\documentclass[11pt,leqno]{article}
\usepackage{graphicx, amsfonts, amsthm, amsxtra, amssymb, verbatim, makeidx}
\usepackage{subeqnarray, relsize}
\usepackage[mathscr]{euscript}
\usepackage{hyperref}
\makeindex
\makeglossary
\begin{document}
\newtheorem{theo}{Theorem}
\newtheorem{exam}{Example}
\newtheorem{coro}{Corollary}
\newtheorem{defi}{Definition}
\newtheorem{prob}{Problem}
\newtheorem{lemm}{Lemma}
\newtheorem{prop}{Proposition}
\newtheorem{rem}{Remark}
\newtheorem{conj}{Conjecture}
\newtheorem{calc}{}

\def\gru{\mu} 
\def\pg{{ \sf S}}               
\def\TS{{\mathlarger{\bf T}}}                
\def\NB{{\mathlarger{\bf N}}}
\def\group{{\sf G}}
\def\NLL{{\rm NL}}   

\def\plc{{ Z_\infty}}    
\def\pola{{u}}      
\newcommand\licy[1]{{\mathbb P}^{#1}} 
\newcommand\aoc[1]{Z^{#1}}     
\def\HL{{\rm Ho}}     
\def\NLL{{\rm NL}}   

\def\Z{\mathbb{Z}}                   
\def\Q{\mathbb{Q}}                   
\def\C{\mathbb{C}}                   
\def\N{\mathbb{N}}                   
\def\uhp{{\mathbb H}}                
\def\A{\mathbb{A}}                   
\def\dR{{\rm dR}}                    
\def\F{{\cal F}}                     
\def\Sp{{\rm Sp}}                    
\def\Gm{\mathbb{G}_m}                 
\def\Ga{\mathbb{G}_a}                 
\def\Tr{{\rm Tr}}                      
\def\tr{{{\mathsf t}{\mathsf r}}}                 
\def\spec{{\rm Spec}}            
\def\ker{{\rm ker}}              
\def\GL{{\rm GL}}                
\def\ker{{\rm ker}}              
\def\coker{{\rm coker}}          
\def\im{{\rm Im}}               
\def\coim{{\rm Coim}}            
\def\p{{\sf  p}}
\def\U{{\cal U}}   

\def\weig{{\nu}}
\def\r{{ r}}                       
\def\k{{\sf k}}                     
\def\ring{{\sf R}}                   
\def\X{{\sf X}}                      
\def\Ua{{   L}}                      
\def\T{{\sf T}}                      
\def\asone{{\sf A}}                  

\def\Ts{{\sf S}}
\def\cmv{{\sf M}}                    
\def\BG{{\sf G}}                       
\def\podu{{\sf pd}}                   
\def\ped{{\sf U}}                    
\def\per{{\bf  P}}                   
\def\gm{{  A}}                    
\def\gma{{\sf  B}}                   
\def\ben{{\sf b}}                    

\def\Rav{{\mathfrak M }}                     
\def\Ram{{\mathfrak C}}                     
\def\Rap{{\mathfrak G}}                     

\def\mov{{\sf  m}}                    
\def\Yuk{{\sf C}}                     
\def\Ra{{\sf R}}                      
\def\hn{{ h}}                         
\def\cpe{{\sf C}}                     
\def\g{{\sf g}}                       
\def\t{{\sf t}}                       
\def\pedo{{\sf  \Pi}}                  

\def\Der{{\rm Der}}                   
\def\MMF{{\sf MF}}                    
\def\codim{{\rm codim}}                
\def\dim{{\rm    dim}}                
\def\Lie{{\rm Lie}}                   
\def\gg{{\mathfrak g}}                

\def\u{{\sf u}}                       

\def\imh{{  \Psi}}                 
\def\imc{{  \Phi }}                  
\def\stab{{\rm Stab }}               
\def\Vec{{\rm Vec}}                 
\def\prim{{\rm  0}}                  

\def\Fg{{\sf F}}     
\def\hol{{\rm hol}}  
\def\non{{\rm non}}  
\def\alg{{\rm alg}}  
\def\tra{{\rm tra}}  

\def\bcov{{\rm \O_\T}}       

\def\leaves{{\cal L}}        

\def\cat{{\cal A}}              
\def\im{{\rm Im}}               

\def\pn{{\sf p}}              
\def\Pic{{\rm Pic}}           
\def\free{{\rm free}}         
\def \NS{{\rm NS}}    
\def\tor{{\rm tor}}
\def\codmod{{\xi}}    

\def\GM{{\rm GM}}

\def\perr{{\sf q}}        
\def\perdo{{\cal K}}   
\def\sfl{{\mathrm F}} 
\def\sp{{\mathbb S}}  

\newcommand\diff[1]{\frac{d #1}{dz}} 
\def\End{{\rm End}}              

\def\sing{{\rm Sing}}            
\def\cha{{\rm char}}             
\def\Gal{{\rm Gal}}              
\def\jacob{{\rm jacob}}          
\def\tjurina{{\rm tjurina}}      
\newcommand\Pn[1]{\mathbb{P}^{#1}}   
\def\P{\mathbb{P}}
\def\Ff{\mathbb{F}}                  

\def\O{{\cal O}}                     

\def\ring{{\mathsf R}}                         
\def\R{\mathbb{R}}                   

\newcommand\ep[1]{e^{\frac{2\pi i}{#1}}}
\newcommand\HH[2]{H^{#2}(#1)}        
\def\Mat{{\rm Mat}}              
\newcommand{\mat}[4]{
     \begin{pmatrix}
            #1 & #2 \\
            #3 & #4
       \end{pmatrix}
    }                                
\newcommand{\matt}[2]{
     \begin{pmatrix}                 
            #1   \\
            #2
       \end{pmatrix}
    }
\def\cl{{\rm cl}}                

\def\hc{{\mathsf H}}                 
\def\Hb{{\cal H}}                    
\def\pese{{\sf P}}                  

\def\PP{\tilde{\cal P}}              
\def\K{{\mathbb K}}                  

\def\M{{\cal M}}
\def\RR{{\cal R}}
\newcommand\Hi[1]{\mathbb{P}^{#1}_\infty}
\def\pt{\mathbb{C}[t]}               
\def\gr{{\rm Gr}}                
\def\Im{{\rm Im}}                
\def\Re{{\rm Re}}                
\def\depth{{\rm depth}}
\newcommand\SL[2]{{\rm SL}(#1, #2)}    
\newcommand\PSL[2]{{\rm PSL}(#1, #2)}  
\def\Resi{{\rm Resi}}              

\def\L{{\cal L}}                     
\def\Aut{{\rm Aut}}              
\def\any{R}                          
\newcommand\ovl[1]{\overline{#1}}    

\newcommand\mf[2]{{M}^{#1}_{#2}}     
\newcommand\mfn[2]{{\tilde M}^{#1}_{#2}}     

\newcommand\bn[2]{\binom{#1}{#2}}    
\def\ja{{\rm j}}                 
\def\Sc{\mathsf{S}}                  
\newcommand\es[1]{g_{#1}}            
\newcommand\V{{\mathsf V}}           
\newcommand\WW{{\mathsf W}}          
\newcommand\Ss{{\cal O}}             
\def\rank{{\rm rank}}                
\def\Dif{{\cal D}}                   
\def\gcd{{\rm gcd}}                  
\def\zedi{{\rm ZD}}                  
\def\BM{{\mathsf H}}                 
\def\plf{{\sf pl}}                             
\def\sgn{{\rm sgn}}                      
\def\diag{{\rm diag}}                   
\def\hodge{{\rm Hodge}}
\def\HF{{ F}}                                
\def\WF{{ W}}                               
\def\HV{{\sf HV}}                                
\def\pol{{\rm pole}}                               
\def\bafi{{\sf r}}
\def\id{{\rm id}}                               
\def\gms{{\sf M}}                           
\def\Iso{{\rm Iso}}                           

\def\hl{{\rm L}}    
\def\imF{{\rm F}}
\def\imG{{\rm G}}

\def\cf{r}   
\def\cm{\checkmark}
\def\MI{{\cal M}}
\def\se{{\sf s}}

\begin{center}
{\LARGE\bf Hodge cycles for cubic hypersurfaces
}
\\
\vspace{.25in} {\large {\sc Hossein Movasati}}
\footnote{
Instituto de Matem\'atica Pura e Aplicada, IMPA, Estrada Dona Castorina, 110, 22460-320, Rio de Janeiro, RJ, Brazil,
{\tt www.impa.br/$\sim$ hossein, hossein@impa.br.}}
\end{center}
\begin{abstract}
We study an algebraic cycle of the form 
$Z_0=\cf\P^{\frac{n}{2}}+\check\cf\check\P^{\frac{n}{2}},\ \cf\in\N,\check\cf\in\Z,\ \  1\leq \cf, |\check\cf|\leq 10,\ \  \gcd (\cf,\check \cf)=1$,
inside the cubic Fermat variety $X_0$ of 
even dimension  $n\geq 4$ and with 
$\dim\left (\P^{\frac{n}{2}}\cap \check\P^{\frac{n}{2}}\right)=m$.
We take a smooth deformation space $\Ts$ of $X_0$ such that the triple  $(X_0, \P^\frac{n}{2}, 
\check\P^\frac{n}{2})$ becomes rigid. 
For $m=\frac{n}{2}-2$ and  for many examples of $N\in\N$ and $n$ 
we show that the $N$-th order Hodge locus attached to  $Z_0$ is smooth and reduced of 
positive dimension if and only if $(\cf,\check\cf)=(1,-1)$. 
In this case, the underlying algebraic cycles are conjectured to be cubic ruled cycles. 
For $m=\frac{n}{2}-3$ the same happens for all choices of coefficients 
$\cf$ and $\check\cf$ and we do not know what kind of algebraic cycles might produce such Hodge cycles. 
The first case gives us a conjectural description of a component of the Hodge locus, and the
second case gives us strong computer assisted evidences for the existence of new Hodge cycles for 
cubic hypersurfaces. 
Whereas the well-known construction of Hodge cycles due to D. Mumford and A. Weil for CM abelian 
varieties, and Y. Andr\'e's motivated cycles can be described  in theoretical terms, 
the  full proof of the existence of our Hodge cycle  seems to be only possible with more  
powerful computing machines.  
\end{abstract}

\section{Introduction}
A. Weil starts his article \cite{Weil1977} with the following:  
`In searching for possible counterexamples to the 
``Hodge conjecture", one has to look for varieties whose Hodge ring is not generated by its elements 
of degree 2'. The description of such elements for CM abelian varieties is fairly understood 
and it is due to D. Mumford and A. Weil himself, see \cite[page 166-167]{Mumford1966-2}, \cite{Weil1977}.
Other examples are Y. Andr\'e's motivated Hodge cycles in \cite{Andre1996}. 
None of these methods can be applied to hypersurfaces. 
By Lefschetz  hyperplane section theorem, for hypersurfaces of even dimension $n\geq 4$ 
we have only a one dimensional subspace of  $H^{\frac{n}{2},\frac{n}{2}}$ generated by elements of 
degree $2$ (the class of a hyperplane section), 
and producing interesting Hodge cycles 
in this case, for which the Hodge conjecture is not known, 
is extremely difficult, and hence, they might be a better candidate for a counterexample to the
Hodge conjecture. Even for Fermat varieties the methods introduced by Z. Ran, T. Shioda in
\cite{Ran1980}, \cite{sh79}  are not based on  explicit description of 
Hodge cycles, and the author in  his book \cite[Chapter 15]{ho13} had to work 
out a computer implementable description of such cycles. This resulted, for instance, to the verification
of integral Hodge conjecture for many examples of the Fermat variety, see \cite{EnzoRoberto}.   
The main goal of Chapter 18 of 
this book is to describe a computer assisted project in order 
to classify components of Hodge loci passing through Fermat point, and in this way to discover new Hodge cycles 
by deforming the class of algebraic cycles. The main difficulty is that 
the parameter space of hypersurfaces are usually of huge dimensions and it is hard to carry out
computations even with the modern computers of today.
\href{http://w3.impa.br/$\sim$hossein/myarticles/Headaches.pdf}{In November 2018, P. Deligne wrote many comments regarding this chapter,}
\footnote{ \href{http://w3.impa.br/$\sim$hossein/myarticles/Headaches.pdf}{See Chapter 1 of 
{\tt w3.impa.br/$\sim$hossein/myarticles/Headaches.pdf.}}}
and this gave the author
some force  to push forward  a tiny step toward one of the main goals of this book. 
This is namely  to find a smaller parameter space
for cubic hypersurfaces and describe instances, where a rigid algebraic
cycle deforms into a  Hodge cycle.

Let $n\geq 4$ be an even number  and  let us consider the projective space  $\P^{n+1}$ with the coordinate system 
$[x_0:x_1:\cdots:x_{n+1}]$ and  $\P^{\frac{n}{2}},\check\P^{\frac{n}{2}}\subset \P^{n+1}$ 
given by:
\begin{equation}
\label{CarolinePilar}
 \P^{\frac{n}{2}}:  
\left\{
 \begin{array}{l}
 x_{0}-\zeta_{6}^{}x_{1}=0,\\
 x_{2}-\zeta_{6}^{} x_{3}=0,\\
 x_{4}-\zeta_{6}^{} x_{5}=0,\\
 \cdots \\
 x_{n}-\zeta_{6}^{} x_{{n+1}}=0,
 \end{array}
 \right.\ \ \ \ \ \ \ 
 \check\P^{\frac{n}{2}}:  
\left\{
 \begin{array}{l}
 x_{0}-\zeta_{6}^{}x_{1}=0,\\
  \cdots \\
 x_{2m}-\zeta_{6}^{} x_{2m+1}=0,\\
 x_{2m+2} + x_{2m+3}=0,\\
 \cdots \\
 x_{n}+ x_{{n+1}}=0,
 \end{array}
 \right.  
 \end{equation}
 where $\zeta_{6}:=e^{\frac{2\pi \sqrt{-1}}{6}}$. 
These are linear algebraic cycles in the cubic Fermat variety $X_0\subset \P^{n+1}$ given  
by the homogeneous polynomial $x_0^3+x_1^3+\cdots+x_{n+1}^3=0$, 
and satisfy $\P^\frac{n}{2}\cap \check\P^\frac{n}{2}=\P^{m}$. For the main purpose of this paper 
we will only consider cubic Fermat varieties and $m=\frac{n}{2}-2, \frac{n}{2}-3$. The other cases are fairly discussed in
\cite[Chapter 18]{ho13}.  We choose a deformation of the Fermat variety
\begin{equation}
\label{15dec2016}
X_t: x_0^{3}+x_1^{3}+\cdots+x_{n+1}^3-\sum_{\alpha \in I}t_\alpha x_{\alpha_1}x_{\alpha_2}x_{\alpha_3}=0,\ \
t=(t_\alpha)_{\alpha\in I}\in\Ts:=\C^{\#I}, 
\end{equation}
where $\alpha$ runs through a finite subset $I$ of all three elements subsets of $\{0,1,\ldots,n+1\}$. 
This deformation is taken in such a way that that the triple $(X_0, \P^\frac{n}{2}, 
\check\P^\frac{n}{2})$ does not deform, see \S\ref{museumofscience2108}. 
For $n=4,6,8,10$  we have computed such a deformation and the 
corresponding monomials are listed in Table \ref{07dec2018-1} and Table \ref{07dec2018-2}. 
For instance, for $n=4$ we have considered the deformations: 
\begin{equation}
\label{07/11/2018-1}
X_t: \ x_0^3+x_1^3+x_2^3+x_3^3+x_4^3+x_5^3-(t_1x_2+t_2x_3)x_1x_5 \hbox{  the case } m=0,
\end{equation}
\begin{equation}
\label{07/11/2018-2}
X_t: \ x_0^3+x_1^3+x_2^3+x_3^3+x_4^3+x_5^3-(t_1x_0+t_2 x_1)x_3x_5,  \hbox{  the case } m=-1.
\end{equation}
For any algebraic cycle $Z_0$ of dimension $\frac{n}{2}$ in $X_0$ let $[Z_0]\in H_n(X_0,\Z)$ be the homology 
class of  $Z_0$ which is a Hodge cycle.  
For a Hodge cycle $\delta_0\in H_n(X_0,\Z)$, and in particular 
$\delta_0$ any linear combination of $[\P^{\frac{n}{2}}]$ and $[\check\P^{\frac{n}{2}}]$
with $\Z$ coefficients,  and $N\in N$, we can define the Hodge locus $V_{\delta_0}$ and $N$-th order 
infinitesimal Hodge locus $V_{\delta_0}^N$,  see \S\ref{10/12/2018}.  The first order Hodge locus $V^1_{\delta_0}$
is the tangent space ${\TS}_0V_{\delta_0}$ of $V_{\delta_0}$ at $0$ and its study is mainly done under the name 
infinitesimal variation of Hodge structures (IVHS) introduced by P. Griffiths and his coauthors in \cite{CGGH1983}.

\begin{theo}
\label{hope2018}
Let  $\P^{\frac{n}{2}}$ and $\check\P^{\frac{n}{2}}$ be linear cycles in \eqref{CarolinePilar} 
with $\P^{\frac{n}{2}}\cap \check\P^{\frac{n}{2}}=\P^m$ and $n=4,6,8,10,12$.
For $m=\frac{n}{2}-2,\ \frac{n}{2}-3$ consider the family of hypersurfaces \eqref{15dec2016} with monomials coming 
from Table \ref{07dec2018-1} and Table \ref{07dec2018-2}, respectively.
\begin{enumerate}
\item{($m=\frac{n}{2}-2$)}
\label{hope-1}
For all $\cf,\check\cf\in\Z,\ 1\leq |\cf|, |\check\cf|\leq 10$,
the infinitesimal Hodge locus $V_{\cf[\P^{\frac{n}{2}}]+\check\cf[\check\P^{\frac{n}{2}}]}^N$
is smooth at $0$ for the cases in Table \ref{16dec2018Tokyo} with $\cm$ mark,  and further   
for $\cf\not=-\check\cf$ it is not smooth at $0$ for the cases in Table \ref{16dec2018Tokyo} with $X$ mark, and so,  
for all these cases the Hodge locus
$V_{\cf[\P^{\frac{n}{2}}]+\check\cf[\check\P^{\frac{n}{2}}]}$ as analytic scheme is either 
non-reduced 
or its underlying analytic variety is singular at the Fermat point 
$0$. Moreover, 
the infinitesimal Hodge locus $V_{[\P^\frac{n}{2}]-[\check\P^\frac{n}{2}]}^N$ 
is smooth for all $N$'s listed in the last row of Table \ref{16dec2018Tokyo}.
\item {($m=\frac{n}{2}-3$)}
\label{hope-2}
For all $\cf,\check\cf\in\Z,\ 1\leq |\cf|, |\check\cf|\leq 10$,
the infinitesimal Hodge locus $V_{\cf[\P^{\frac{n}{2}}]+\check\cf[\check\P^{\frac{n}{2}}]}^N$
is  smooth at $0$ for all $N$'s listed in Table \ref{nabeta2018}. Moreover, their tangent spaces 
form a pencil with the origin as axis, that is, they have the same dimension and intersect
each other only at the origin. 
\item
\label{hope-3}
The codimension of the Zariski tangent space of $V_{\cf[\P^\frac{n}{2}]+\check\cf[\check\P^\frac{n}{2}]}$ at $0$ in both 
cases as above is listed in the third row  of  Table \ref{16dec2018Tokyo} and  Table \ref{nabeta2018}, respectively. 
\item
\label{hope-4}
The Fermat variety $X_0$ together with its algebraic cycles $\P^\frac{n}{2},\ \check\P^\frac{n}{2}$ with
$m=\frac{n}{2}-2$ (resp. $m=\frac{n}{2}-3$) is 
rigid inside the family \eqref{15dec2016} with monomials coming  from Table \ref{07dec2018-1} 
(resp. \ref{07dec2018-2}.
\end{enumerate}
\end{theo}

\begin{table}[!htbp]
\centering
\begin{tabular}{|c|c|c|c|c|c|}
 \hline  
\multicolumn{1}{|c|}{ $n$} & $4$   & $6$    & $8$    & $10$  & $12$ \\ \hline\hline
\multicolumn{1}{|c|}{$\dim(\Ts)$}      &  $2$  & $8$     &  $19$  & $36$  & $60$   \\ \hline                                                   
\multicolumn{1}{|c|}{$\codim( V_{\cf[\P^\frac{n}{2}]+\check\cf[\check\P^\frac{n}{2}]}^1      )$}&$1$ & $6$     & $16$   & $32$  & $55$  \\ \hline\hline 
\multicolumn{1}{|c|}{$N=2$}             & $\cm$ & $\cm$   & $\cm$  & $\cm$ &    $\cm$   \\ \hline
\multicolumn{1}{|c|}{$N=3$}             & $\cm$ & $\cm$   & $X$  & $X$     &    $X$   \\ \hline
\multicolumn{1}{|c|}{$N=4$}             & $\cm$ & $X$   & $X$  & $X$     &   $X$    \\ \hline
\hline
\multicolumn{1}{|c|}{$Z_0=\P^{\frac{n}{2}}-\check\P^{\frac{n}{2}},\ \  \ N=$}
                                      &
                                        \href{http://w3.impa.br/~hossein/WikiHossein/files/Singular%20Codes/2018-12_CubicHypersurfaces_n=4_d=3_m=0_tru=20.txt}{$\infty$}
                                              & 
                                               \href{http://w3.impa.br/~hossein/WikiHossein/files/Singular%20Codes/2018-12_CubicHypersurfaces_n=6_d=3_m=1_tru=7.txt}{$14$} 
 
                                                        & 
                                                         \href{http://w3.impa.br/~hossein/WikiHossein/files/Singular%20Codes/2018-12_CubicHypersurfaces_n=8_d=3_m=2_tru=4.txt}{$6$}
                                                                  & $4$  &   $3$      \\ \hline

\end{tabular}
\label{16dec2018Tokyo}
\caption{Smooth and singular Hodge loci: $m=\frac{n}{2}-2$}
\end{table}

\begin{table}[!htbp]
\centering
\begin{tabular}{|c|c|c|c|c|c|}
 \hline  
\multicolumn{1}{|c|}{ $n$}                    & $4$   &  $6$ & $8$  &  $10$ &  $12$ \\ \hline \hline
\multicolumn{1}{|c|}{$\dim(\Ts)$}              & $2$   &$8$   & $20$ & $39$  & $66$      \\ \hline                                                   
\multicolumn{1}{|c|}{$\codim( V_{\cf[\P^\frac{n}{2}]+\check\cf[\check\P^\frac{n}{2}]}^1)$  }   &$1$    &$7$   & $19$  & $38$  & $65$      \\ \hline\hline 
\multicolumn{1}{|c|}{$Z_0=\cf\P^{\frac{n}{2}}+\check\cf\check\P^{\frac{n}{2}},\ \ N=$}
                                              &
                                                                                 \href{  http://w3.impa.br/~hossein/WikiHossein/files/Singular%20Codes/2018-12_CubicHypersurfaces_n=4_d=3_m=-1_tru=20.txt  }{$\infty$}
                                                      & 
                                                                                        \href{http://w3.impa.br/~hossein/WikiHossein/files/Singular%20Codes/2018-12_CubicHypersurfaces_n=6_d=3_m=0_tru=7.txt}{$14$}

                                                             &  
                                                                                      \href{http://w3.impa.br/~hossein/WikiHossein/files/Singular%20Codes/2018-12_CubicHypersurfaces_n=8_d=3_m=1_tru=4.txt}{$6$}
 
                                                                    &  
                                                                          $4$    &      $3$
                                                                                                                    \\ \hline

\end{tabular}
\caption{Smooth  Hodge loci: $m=\frac{n}{2}-3$}
\label{nabeta2018}
\end{table}

Theorem \ref{hope2018} and similar computation in \cite[Chapter 18]{ho13} for the full family of hypersurfaces lead us to 
speculate many conjectures. The following might be the most
evident one. 
\begin{conj}
\label{doran2018}
\begin{enumerate}
 \item 
For the full family of cubic hypersurfaces the Hodge locus 
$V_{[\P^\frac{n}{2}]-[\check\P^\frac{n}{2}]}$ with $\P^\frac{n}{2}\cap \check\P^\frac{n}{2}=
\P^{\frac{n}{2}-2}$ 
is smooth and it is larger than the deformation space of the triple 
$(X_0,\P^\frac{n}{2}, \ \check\P^\frac{n}{2})$. 
\item
For the full family of cubic hypersurfaces the Hodge locus 
$V_{\cf[\P^\frac{n}{2}]-\check\cf[\check\P^\frac{n}{2}]}$ with $\P^\frac{n}{2}\cap \check\P^\frac{n}{2}=
\P^{\frac{n}{2}-3}$ and $|\cf|,|\check\cf|\in\N$,  
is smooth and it is larger then the deformation space of the triple 
$(X_0,\P^\frac{n}{2}, \ \check\P^\frac{n}{2})$. The difference of dimensions in this case is $1$. 
\end{enumerate}
\end{conj}
Conjecture \ref{doran2018} is true for cubic fourfolds for trivial reasons.
In this case $h^{40}=0,\ h^{31}=1$ and 
so any Hodge locus $V_{\delta_0}$, $\delta_0\in H_n(X_0,\Z)$ primitive non-zero Hodge cycle, 
is given by one equation, which turns out that it has non-zero linear 
part when we consider the full parameter space of cubic hypersurfaces, and hence it is always smooth.  Both the integral and rational Hodge conjecture are proved in this case, 
see \cite{zu77} and \cite[Theorem 2.11 and the comments thereafter]{voisin2013}.
In this case, an effective construction of algebraic cycles for hypersurfaces parameterized by 
$V_{\cf [\P^2]+\check\cf[\check\P^2]}$ might give some hint to do a similar verification
in higher dimensions. For a partial verification of Conjecture \ref{doran2018} part 1 see \S\ref{30jan2019}.


\begin{table}
\centering
\begin{tabular}{|c|c|}
\hline   
 $n=4$       & $x_1x_2x_5$, $x_1x_3x_5$   \\  \hline
 $n=6$       &  $x_{1}x_{3}x_{4}$, $x_{1}x_{3}x_{5}$, 
                $x_{1}x_{3}x_{6}$, $x_{1}x_{3}x_{7}$, 
                $x_{1}x_{4}x_{7}$, $x_{3}x_{4}x_{7}$, 
                $x_{1}x_{5}x_{7}$, $x_{3}x_{5}x_{7}$  
                \\  \hline
 $n=8$       &    
 $x_{1}x_{3}x_{5}$
$x_{1}x_{3}x_{6}$
$ x_{1}x_{5}x_{6}$
$ x_{3}x_{5}x_{6}$
$ x_{1}x_{3}x_{7}$
$ x_{1}x_{5}x_{7}$
$ x_{3}x_{5}x_{7}$
$ x_{1}x_{3}x_{8}$
$ x_{1}x_{5}x_{8}$  \\ 
     &                  
$ x_{3}x_{5}x_{8}$
$ x_{1}x_{3}x_{9}$
$ x_{1}x_{5}x_{9}$
$ x_{3}x_{5}x_{9}$
$ x_{1}x_{6}x_{9}$
$ x_{3}x_{6}x_{9}$
$ x_{5}x_{6}x_{9}$
$ x_{1}x_{7}x_{9}$
$ x_{3}x_{7}x_{9}$
$ x_{5}x_{7}x_{9}$
                 \\  \hline
 $n=10$        &    
$x_{1}x_{3}x_{5}$, 
$x_{1}x_{3}x_{7}$,
$x_{1}x_{5}x_{7}$
$x_{3}x_{5}x_{7}$
$x_{1}x_{3}x_{8}$
$x_{1}x_{5}x_{8}$
$x_{3}x_{5}x_{8}$
$x_{1}x_{7}x_{8}$
$x_{3}x_{7}x_{8}$  \\ 
  &                
$x_{5}x_{7}x_{8}$
$x_{1}x_{3}x_{9}$
$ x_{1}x_{5}x_{9}$
$x_{3}x_{5}x_{9}$
$x_{1}x_{7}x_{9}$
$x_{3}x_{7}x_{9}$
$x_{5}x_{7}x_{9}$
$x_{1}x_{3}x_{10}$
$x_{1}x_{5}x_{10}$ \\ 
   & 
$x_{3}x_{5}x_{10}$
$x_{1}x_{7}x_{10}$
$x_{3}x_{7}x_{10}$
$x_{5}x_{7}x_{10}$
$x_{1}x_{3}x_{11}$
$x_{1}x_{5}x_{11}$
$x_{3}x_{5}x_{11}$
$x_{1}x_{7}x_{11}$
$x_{3}x_{7}x_{11}$  \\
   &
$x_{5}x_{7}x_{11}$ 
$x_{1}x_{8}x_{11}$
$x_{3}x_{8}x_{11}$
$x_{5}x_{8}x_{11}$
$x_{7}x_{8}x_{11}$
$x_{1}x_{9}x_{11}$
$x_{3}x_{9}x_{11}$
$x_{5}x_{9}x_{11}$
$ x_{7}x_{9}x_{11}$
                   \\  \hline
$n=12$ &  
\href{http://w3.impa.br/~hossein/WikiHossein/files/Singular%20Codes/2018-12-60monomials-DeformationOfCubicDegree12.txt}
{$60$ monomials}  \\  \hline                    
\end{tabular}
\caption{Monomials of a deformation: $m=\frac{n}{2}-2$} 
\label{07dec2018-1}
\begin{center}
 
\end{center}

\begin{tabular}{|c|c|}
\hline   
 $n=4$       &  $x_0x_3x_5$, $x_1x_3x_5$  \\  \hline
 $n=6$       &   
              $x_1x_2x_5$, $x_1x_3x_5$, $x_1x_2x_7$, $x_1x_3x_7$, $x_1x_4x_7$, 
              $x_1x_5x_7$, $x_2x_5x_7$, $x_3x_5x_7$
  \\  \hline
 $n=8$       &   $x_1x_3x_4$, $x_1x_3x_5$, $x_1x_3x_6$, $x_1x_3x_7$, $x_1x_4x_7$, $x_3x_4x_7$, $x_1x_5x_7$, $x_3x_5x_7$, $x_1x_3x_8$ \\
            &    $x_1x_3x_9$, $x_1x_4x_9$, $x_3x_4x_9$, $x_1x_5x_9$, $x_3x_5x_9$, $x_1x_6x_9$, $x_3x_6x_9$, $x_1x_7x_9$, $x_3x_7x_9$ \\
            &    $x_4x_7x_9$, $x_5x_7x_9$
                 \\  \hline
 $n=10$        &    
                 $x_1x_3x_5$, $x_1x_3x_6$, $x_1x_5x_6$, $x_3x_5x_6$, $x_1x_3x_7$, $x_1x_5x_7$, $x_3x_5x_7$, $x_1x_3x_8$, $x_1x_5x_8$ \\
               & $x_3x_5x_8$, $x_1x_3x_9$, $x_1x_5x_9$, $x_3x_5x_9$, $x_1x_6x_9$, $x_3x_6x_9$, $x_5x_6x_9$, $x_1x_7x_9$, $x_3x_7x_9$ \\
               & $x_5x_7x_9$, $x_1x_3x_{10}$, $x_1x_5x_{10}$, $x_3x_5x_{10}$, $x_1x_3x_{11}$, $x_1x_5x_{11}$, $x_3x_5x_{11}$, $x_1x_6x_{11}$, $x_3x_6x_{11}$ \\
               & $x_5x_6x_{11}$, $x_1x_7x_{11}$, $x_3x_7x_{11}$, $x_5x_7x_{11}$, $x_1x_8x_{11}$, $x_3x_8x_{11}$, $x_5x_8x_{11}$, $x_1x_9x_{11}$, $x_3x_9x_{11}$ \\
               & $x_5x_9x_{11}$, $x_6x_9x_{11}$, $x_7x_9x_{11}$
                   \\  \hline
$n=12$ &  
\href{http://w3.impa.br/~hossein/WikiHossein/files/Singular%20Codes/2018-12-66monomials-DeformationOfCubicDegree12.txt}
{$66$ monomials}  \\  \hline                    
\end{tabular}
\caption{Monomials of a deformation for $m=\frac{n}{2}-3$} 
\label{07dec2018-2}
\end{table}
For the proof of Theorem \ref{hope2018}, the author has written many procedures which are collected in the library  
{\tt foliation.lib}
of {\sc Singular},  see \cite{GPS01}. In order to check the computations of the present paper,  we 
first get this library  from \href{http://w3.impa.br/~hossein/mod-fol-exa.html}
{the author's web page. }
\footnote{\tt http://w3.impa.br/$\sim$hossein/foliation-allversions/foliation.lib}
Then we run the example session of the procedure. For instance, for the procedure {\tt InterTang} used in 
\S\ref{museumofscience2108} we run
\begin{verbatim}
  LIB foliation.lib;
  example InterTang; 
\end{verbatim}
Modifying, the code in the example session (for instance changing the dimension $n$ or degree $d$ of the hypersurface)
we  get all the claimed statements.  

The present work would not have been possible without the attention of two great mathematicians: thanks go
S.-T. Yau for all his effort to create lovely ambients to do mathematics, from CMSA to TSIMF which I enjoyed both institutes 
during the preparation of this text,
and to P. Deligne for all his enlightening emails and comments to the author's book and this article. I would also like
to thank D. van Straten for his help in \S\ref{30jan2019}.

\section{Smaller deformation space}
\label{museumofscience2108}
Let $\T$ be the (full) parameter space of smooth 
cubic hypersurfaces of dimension $n$ and $\X/\T$ be the corresponding family. 
Let also $X_0, \ 0\in\T$ be a smooth hypersurface given by the zero set of a homogeneous 
polynomial of the form
 \begin{equation}
\label{19nov2015-1}
f=f_1f_{s+1}+f_2f_{s+2}+\cdots+f_sf_{2s},\ \ 
f_i\in \C[x]_{d_i},\ f_{s+i}\in  \C[x]_{d-d_i},\ \ s:=\frac{n}{2}+1.
\end{equation}
We call $Z_0: f_1=f_2=\cdots=f_s=0$ a complete intersection algebraic cycle. 
Let $V_{Z_0}\subset (\T,0)$ be the analytic variety  parameterizing  deformations of $(X_0,Z_0)$.
This corresponds to variation of polynomials $f_i$ as above.
This is a branch of the algebraic set $\T_{\underline {d}}\subset \T$ which parameterizes all 
hypersurfaces given by $f$ of the form \ref{19nov2015-1}. For $d_1=d_2=\cdots=d_s=1$, 
$\T_{\underline{1}}$ has $N:=1\cdot 3\cdots (n-1)(n+1) d^{\frac{n}{2}+1}$ branches near the Fermat point 
$0\in\T$, see \cite[\S 17.4]{ho13}. If $f_1,f_2,\cdots, f_{2s}$ have not common 
zeros in $\P^{n+1}$ then they form a regular sequence and a Koszul complex argument tells us that 
$V_Z$ is smooth at $0$ and its tangent space at this point is given by the degree $d$ part of 
the homogeneous ideal
$\langle f_1,f_2,\ldots,f_{2s}\rangle$, see \cite[Proposition 17.5]{ho13}. In particular, this
is the case for two linear cycles 
\begin{eqnarray*}
\P^{\frac{n}{2}} &:&  f_1=f_2=\cdots=f_{s}=0,\\
\check\P^{\frac{n}{2}} &:& \check f_1=\check f_2=\cdots=\check f_{s}=0,\\
\end{eqnarray*}
inside the Fermat variety $X_0$ and given in \eqref{CarolinePilar}.
Moreover, the intersection $V_{\P^{\frac{n}{2}}+ \check\P^{\frac{n}{2}} }:=
V_{\P^{\frac{n}{2}} }\cap V_{\check\P^{\frac{n}{2}} }$ is smooth at $0$ and its 
tangent space at this point is given by the intersection of tangent spaces of 
$V_{\P^{\frac{n}{2}} }$ and $V_{\check\P^{\frac{n}{2}}}$ at $0$,  
\cite[Proposition 17.8]{ho13}. We conclude that 
$$
\TS_0V_{\P^{\frac{n}{2}}+ \check\P^{\frac{n}{2}}}=I_d,\ \ 
{\cal I}:=\langle f_1,f_2,\ldots,f_{2s}\rangle\cap \langle \check f_1,\check f_2,\ldots,\check f_{2s}\rangle.
$$
We make the identification of $\T$ with a Zariski open subset of $\C[x]$ such that $f\in \T$ parametrizes 
the hypersurface given by $x_0^d+x_1^d+\cdots+x^d_{n+1}+f=0$. In this way, $\TS_0\T=\C[x]$.  
We choose a monomial basis $x^\alpha,\ \alpha\in I$ of  degree $d$ piece of the quotient $\C[x]/{\cal I}$ 
and   define $\Ts$ to be a linear subspace of $\C[x]$ generated by these monomials. By definition 
it is perpendicular to ${\cal I}_d$, and this is all what we need. 
The monomials in Table \ref{07dec2018-1} and Table \ref{07dec2018-2} are obtained
in this way.  
This has been implemented in the procedure {\tt InterTang}.

\section{Infinitesimal Hodge loci}
\label{10/12/2018}
After P. Griffiths' work \cite{gr69} we know that for cubic hypersurfaces $X_t$ with  $t$ in a Zariski 
neighborhood of $0\in\T$, the primitive de Rham cohomology $H^{n}_\dR(X_t)_\prim$ has a basis of the form 
\begin{eqnarray*}
& & \omega_\beta := 
\Resi\left(\frac{x_{\beta_1}x_{\beta_2}\cdots x_{\beta_{3k-n-2}}\Omega }{f^k}\right), \hbox{ where }\\
& &    
\Omega=\sum_{i=0}^{n+1}(-1)^ix_i dx_0\wedge dx_1\wedge\cdots \wedge dx_{i-1}\wedge dx_i\wedge \cdots \wedge dx_{n+1}, 
\\
& & 
k\in\N,\ \ \frac{n+2}{3}\leq k\leq \frac{2(n+2)}{3},  \ \ 
\beta:=\{\beta_1, \beta_2,\ldots, \beta_{3k-n-2}\}\subset \{0,1,2,\ldots,n+1\}, \\
\end{eqnarray*}
where $\Resi : H^{n+1}_\dR(\Pn {n+1}-X_t)\to H^n_\dR(X_t)_\prim$  is the residue map, $\beta_i$'s are distinct and the subindex $\prim$ refers
to primitive cohomology.
Moreover, this basis is compatible with the Hodge filtration, that is, 
$F^{n+1-a}H^{n}_\dR(X_t)_\prim$ is generated by $\omega_{\beta}$ with $k\leq a$. 
For the definition of a Hodge locus we need $a=\frac{n}{2}$.
For a 
Hodge cycle $\delta_0\in H_{n}(X_0,\Z)$, the Hodge locus $V_{\delta_0}\subset(\T,0)$ is 
an analytic scheme given by the ideal 
\begin{equation}
\label{17dec2018}
\left\langle\ \ \  \mathlarger{\int}_{\delta_t}\omega_\beta \ \ \ \  \Bigg|   
\ \ \ \forall \beta \subset \{0,1,2,\ldots,n+1\} \hbox{ with },\ \  \#\beta=3k-n-2, \ \  k\leq \frac{n}{2}
\right\rangle\subset \O_{\T,0}.
\end{equation}
Let $\MI_{\T,0}$ be the maximal ideal of $\O_{\T,0}$, that is, the set of germs of holomorphic functions in $(\T,0)$
vanishing at $0$. The $N$-th order infinitesimal scheme  
$V^N_{\delta_0}$ is the induced scheme  by \eqref{17dec2018} in the infinitesimal scheme $\T^N:=
\spec(\O_{\T,0}/\MI^{N+1}_{\T,0})$. We denote by $\X^N/\T^N$ the $N$-th order 
infinitesimal deformation of $X_0$ induced by $\X/\T$. 
Let $\cl(Z_0)\in H^n_\dR(X_0)$ be the class of an algebraic
cycle $Z_0$ of codimension $\frac{n}{2}$ in $X_0$. Let us consider the  Gauss-Manin connection 
$$
\nabla: H^n_\dR(\X/\T)\to \Omega_\T^1\otimes_{\O_\T} 
H^m_\dR(\X/\T).
$$
It induces a connection in $H^n_\dR(\X^N/\T^N)$ which we call it again the Gauss-Manin connection.  
There is a unique section $\se$ 
of  $H^n_\dR(\X^N/\T^N)$ such that $\nabla(\se)=0$ and $\se_0=\cl(Z_0)$. This is called
the horizontal extension of $\cl(Z_0)$ or a flat section of the cohomology bundle.
An equivalent definition for $V^N_{[Z_0]}$ is as follows.
\begin{defi}\rm
\label{kisin2018} 
 The  Hodge locus $V_{[Z_0]}^N$ is a subscheme of $\T^N$ 
 given by the conditions
\begin{eqnarray}
 \label{17/11/2018-1}
& & \nabla(\se)=0,\\ \label{17/11/2018-2}
&  & \se \in F^{\frac{n}{2}}H^n_\dR(\X^N/\T^N),\\  \label{17/11/2018-3}
& & \se_0=\cl(Z_0).
\end{eqnarray}
\end{defi}

\section{Proof of Theorem \ref{hope2018}}
The main ingredient of the proof is a closed formula for the Taylor expansion of the holomorphic
functions $\int_{\delta_t}\omega_\beta$. Since we have only used  a computer implementation of this formula, we do not reproduce it
here and refer the reader to \cite[\S 18.5, Chapter 19]{ho13}. The proof of Item 1 and Item 2 are the same
as the proof of \cite[Theorem 18.2, Theorem 18.3]{ho13}. 
For Item \ref{hope-3} note that these numbers are the dimension of the $\C$-vector space generated by
linear parts of $\int_{\delta_t}\omega_\beta$ in \eqref{17dec2018}.  
Item \ref{hope-4} is just the consequence of our construction of $\Ts$ in 
\S\ref{museumofscience2108}. Note that in a neighborhood of $0$ in $\Ts$, $\Ts$ intersects 
$V_{\P^{\frac{n}{2}}+\check \P^{\frac{n}{2}}}$ only at $0$. 
The following computer code has been used for the proof of statements in 
Theorem \ref{hope2018} Item 1 and Item 2 involving arbitrary coefficients. A simple
modification of it can be used for the case $(\cf,\check\cf)=(1,-1)$.
{\tiny 
\begin{verbatim}
  LIB "foliation.lib"; 
  int n=10; int d=3; int m=(n div 2)-2; int tru=4; int zb=10; 
  intvec zarib1=1,-zb; intvec zarib2=zb, zb;
  intvec mlist=d;  for (int i=1;i<=n; i=i+1){mlist=mlist,d;}
  ring r=(0,z), (x(1..n+1)),dp;
  poly cp=cyclotomic(2*d); int degext=deg(cp) div deg(var(1));
  cp=subst(cp, x(1),z); minpoly =number(cp);
  list lcycles=SumTwoLinearCycle(n,d,m,1); lcycles;
  list ll=InterTang(n,d,lcycles);  
  SmoothReduced(mlist,tru, lcycles, zarib1, zarib2,lcycles);
\end{verbatim}
}

%
%
%

The tangent spaces of $V_{[\P^{\frac{n}{2}}]}= V_{\P^{\frac{n}{2}}}$ and 
$V_{[\check\P^{\frac{n}{2}}]}=V_{\check\P^{\frac{n}{2}}}$ are of the form
$\ker(A)$ and $\ker(\check A)$, where $A$ and $\check A$ are two 
$\dim(\Ts)\times h^{\frac{n}{2}+1,\frac{n}{2}-1}$
matrices with entries in $\Q(\zeta_6)$ which can be constructed from the periods of the underlying 
algebraic cycles, see \cite[\S 16.5]{ho13}. The tangent space of
$V_{[\P^{\frac{n}{2}}]+x [\check\P^{\frac{n}{2}}]},\ x\in\Q$ at $0$ is given by $\ker(A+x \check A)$.
By our construction of $\Ts$, $\ker(A)\cap \ker(\check A)=\{0\}$ which is the tangent space
of the deformation space of $(X_0,\P^{\frac{n}{2}}, \check \P^{\frac{n}{2}} )$. For two distinct 
$x_1,x_2\in\Q$, we have $\ker(A+x_1\check A)\cap \ker(A+x_2\check A)=
\ker(A)\cap \ker(\check A)=\{0\}$ and the last affirmation in Item \ref{hope-2} follows.
Note that in all the cases considered in Table \ref{nabeta2018}, $\ker(A+x_1\check A)$ is a one 
dimensional space which is a byproduct of our computations above.

\section{Finding algebraic cycles}
\label{30jan2019}
After the first draft of this paper was written, there was many email exchanges in January 2019 with P. Deligne in order
to verify the Hodge conjecture for the Hodge cycles in Theorem \ref{hope2018} part \ref{hope-1}. 
Most of the content of this section is the result of this joint effort. We were not able to give similar 
descriptions for Theorem \ref{hope2018} part \ref{hope-2}. Let 
\begin{equation}
 \label{CarolinePilar-1}
 \P^{\frac{n}{2}}_{a_1, a_2}:  
\left\{
 \begin{array}{l}
 x_{0}-\zeta_{6} x_{1}=0,\\
  \cdots \\
 x_{n-4}-\zeta_{6} x_{n-3}=0,\\
 x_{n-2} -\zeta_6^{2a_1+1}  x_{n-3}=0,\\
 x_{n} -\zeta_6^{2a_2+1}  x_{{n+1}}=0.
 \end{array}
 \right.  \hbox{ where }  0\leq a_1,a_2\leq 2.
\end{equation}
Using this notation we have  
$\P^{\frac{n}{2}}_{0, 0}= \P^{\frac{n}{2}}, \P^{\frac{n}{2}}_{1,1}= \check\P^{\frac{n}{2}}$. 
It turns out that 
\begin{equation}
\label{detroitairport2019}
\P^{\frac{n}{2}}_{0, 0}-\P^{\frac{n}{2}}_{1, 1}=
\P^{\frac{n}{2}}_{0,0}+\P^{\frac{n}{2}}_{0,1}+\P^{\frac{n}{2}}_{2,1}
\end{equation}
which is written modulo $\P^{\frac{n}{2}+1}$ slices of the Fermat $X^3_n$, and hence,
both sides of  \eqref{detroitairport2019} induce the same element in the primitive cohomology and have
the same Hodge locus. Let $C_0$ be the algebraic cycle in the right hand side of 
\eqref{detroitairport2019}. 
We write the three cycles in $C_0$ as
\begin{eqnarray*}
 \P^\frac{n}{2}_{0,0} &:& g_1=g_2=\cdots=g_{\frac{n}{2}-1}=f_{11}=f_{21}=0, \\
 \P^\frac{n}{2}_{0,1} &:& g_1=g_2=\cdots=g_{\frac{n}{2}-1}= f_{11}=f_{32}=0,  \\
 \P^\frac{n}{2}_{2,1} &:& g_1=g_2=\cdots=g_{\frac{n}{2}-1}=f_{22}=f_{32}=0,
\end{eqnarray*}
where $g_i$ and $f_{ij}$'s are homogeneous linear polynomials. Now, we can easily see that this algebraic 
cycle deforms into:
\begin{equation}
\label{shodam2019}
C: g_1=g_2=\cdots=g_{\frac{n}{2}-1}=0, \ \ \ 
\rank
\begin{bmatrix}
 f_{11} & f_{12}\\
  f_{21} & f_{22}\\
 f_{31} & f_{32}\\
\end{bmatrix}
\leq 1,
\end{equation}
where $g_i$'s and $f_{ij}$'s are deformed homogeneous linear polynomials, and for simplicity,
we have not introduced new notation for deformed polynomials. 
The cycle $C$ deforms to $C_0$ by setting $f_{12},f_{31}$ equal to zero.
In a personal communication 
Duco van Straten pointed out the determinantal structure of $C$ and the fact that  
for $\frac{n}{2}=2$, $C$ is
the cubic ruled surface/Hirzebruch surface $F_1$. It is isomorphic to $\P^2$
blown up in a single point, embedded by the linear
system of quadrics through the point. For this reason, we call $C$ a cubic ruled cycle of dimension 
$\frac{n}{2}$. 

The ideal given by the determinantal variety in \eqref{shodam2019} is radical, and hence, if a 
smooth hypersurface
of degree $d$ and dimension $n$ given by the homogeneous polynomial $f$ contains $C$ then 
\begin{equation}
\label{21jan2019}
f=g_1*_1+g_2*_2+\cdots+g_{\frac{n}{2}-1}*_{\frac{n}{2}-1}+\left| \begin{matrix}
     f_{11} & f_{12} & f_{13}\\
     f_{21} & f_{22} & f_{23}\\
     f_{31} & f_{32} & f_{33}
    \end{matrix}
\right|,
\end{equation}
where $*_i$'s are homogeneous polynomials of degree $d-\deg(g_i)$ and the first two columns consist 
of linear homogeneous polynomials and the last column of degree $d-2$ 
homogeneous polynomials. Let $\check\T\subset\T$ be the space of smooth hypersurfaces of the form 
\eqref{21jan2019} and $V_{C_0}\subset \check\T$ be the branch of $\check\T$ corresponding to deformations
of $(X_0,C_0)$. It is now natural to conjecture that $V_{[C_0]}=V_{C_0}$. For this it is enough 
to prove that the codimension of $\check \T$ in $\T$ is the same as the codimension of 
$\TS_0V_{[C_0]}$  in $\TS_0\T$ which is computed  in Table \ref{16dec2018Tokyo}. 
The author was not able to find a 
closed  formula for the codimension of $\check\T$. Instead, it is possible 
to take random points in $\check\T$ and compute the tangent space of $\check\T$
at such points for all examples of $n$ in Table \ref{16dec2018Tokyo}.  
It turned out that one gets the desired codimension.  For this computer assisted 
verification, we have used {\tt CodRuledCubic}.  

It might be possible to prove that $\check\T$ is a component of the Hodge locus, that is, for generic 
$X_t,\ t\in\check\T$ with $C\subset X$ as above, the homology class of $C$ cannot be deformed into a 
Hodge cycle in $X_t, t\in \T\backslash\check\T$. For sum of two linear cycles such statements 
are proved  in \cite{Roberto} using computer calculations and in \cite{RobertoThesis} using long  
theoretical methods.

It is early to claim that the conjectural Hodge cycles introduced in Theorem 
\ref{hope2018} part \ref{hope-2} are counterexamples to the Hodge conjecture, as this needs 
more time and effort of other people who might try to verify the Hodge conjecture for
such cycles. However, it seems to the author that they are
better candidates for this than the well-known ones. 
Components of low codimension of the Noether-Lefschetz loci (Hodge loci for $n=2$) parameterize surfaces
with rather  simple algebraic cycles such as lines or quadrics, see for instance 
\cite{green1989, voisin1988}.  The minimal codimension of components of the Hodge loci 
for cubic hypersurfaces in a Zariski neighborhood of the
Fermat point $0\in\T$ is 
$\binom{\frac{n}{2}+3}{3}-(\frac{n}{2}+1)^2=\binom{\frac{n}{2}+1}{3}$ which is obtained 
by $\T_{\underline{1}}$ introduced in 
\S\ref{museumofscience2108},  see \cite[Theorem 2]{GMCD-NL}. 
The next admissible codimension seems to be of $\check\T$ introduced in this section. 
Note that for cubic hypersurfaces all the components $\T_{\underline{d}}$ are
the same as $\T_{\underline{1}}$. After this, we have conjecturally an infinite number of   
components $V_{\cf\P^{\frac{n}{2}}+\check\cf\P^{\frac{n}{2}} },\ \ m=\frac{n}{2}-3$ of the 
same codimension. 
For the convenience of the reader we have also computed Table \ref{12mar2017-2} which contains 
the dimension of the full moduli,  Hodge numbers and  
the range of codimensions of  Hodge loci  for cubic  hypersurfaces. 
Finally, note  that it follows from  \cite[Theorem 2.12, Principle B]{dmos} that
the Hodge cycles of the present text are absolute.

\def\Has{{\cal C}}
\section{Learning from cubic fourfolds}
Despite the fact that the integral Hodge conjecture for cubic fourfolds is well-known, 
effective construction of algebraic cycles in this case might help us in a better understanding
of the Hodge cycles in Theorem \ref{hope2018}. For cubic fourfolds, Hodge loci is a union of  
codimension one irreducible subvarieties $\Has_D,\ D\equiv_6 0,2, D\geq 8$ of $\T$, see 
\cite{Hassett2000}. Here, $D$ is the discriminant of the saturated lattice generated 
by $[Z]$ and $[Z_\infty]$ in $H_4(X,\Z)$ (in \cite{Hassett2000} notation $[Z_\infty]=h^2$), where   
$Z$ is an  algebraic cycle $Z\subset X, \ X\in \Has_D$ whose homology class together $[Z_\infty]$ form a rank two lattice. The loci of cubic fourfolds containing a plane (resp. cubic ruled surface) is 
$\T_{1,1,1}=\Has_8$ (resp. $\Has_{12}$). A cubic ruled surface is also called
cubic scroll. Let $Z=\cf \P^2+\check\cf\check\P^2$ be as in Theorem \ref{hope2018} with $\P^2\cap\check \P^2=\P^m,\ \  m=-1,0,1$, 
$\cf,\check\cf$ are two non-zero coprime integers and $\cf>0$.  Using adjunction formula, see 
for instance \cite[Section 17.6]{ho13}, we can easily see that the discriminant of the lattice generated
by $Z$ and $Z_\infty$ is 
$$
D:=
\left\{
 \begin{array}{ll}
8(\cf^2+\check \cf^2)-2\cf\cdot \check\cf & m=-1,\ \ \hbox{ two disjoint planes,}\\
8(\cf^2+\check \cf^2)+4\cf\cdot \check\cf &  m=0, \ \  \hbox{two planes intersecting in a point,}\\
8(\cf^2+\check \cf^2)-8\cf\cdot \check\cf &  m=1, \ \  \hbox{two planes intersecting in a line.}
\end{array}\right. 
$$
The case $m=-1$ is our main interest and it is  always $\equiv_6 0,2$, however, the corresponding lattice
might not be saturated. 
The first values of $D$ above  are $14,18,36$ which are obtained by 
$(\cf,\check\cf)=(1,1),(1,-1),(2,1)$, respectively. Below, we will describe generalizations of 
$\Has_{14}$ and $\Has_{20}$ for cubic $n$-folds. They seems to be far from any possible algebraic
cycle for the verification of Theorem \ref{hope2018} part \ref{hope-2}. 

The loci $\Has_{14}$ parametrizes cubic fourfolds with a quartic scroll. This is the image of 
$$
\P^1\times \P^1\hookrightarrow \P^5,\ (x,y)\mapsto [f_1g_1:f_1g_2:f_1g_3: f_2g_1:f_2g_2:f_2g_3],
$$ 
where $f_i$'s (resp. $g_i$'s) form a basis of $\C[x,y]_1$ (resp. $\C[x,y]_2$). Putting 
$f_1=x,\ f_2=y,\ g_1=x^2,\ g_2=xy,\ g_3=y^2$, the ideal of a quartic scroll is given by
\begin{equation}
\label{shodam2019}
C:\  f_{21}f_{22}-f_{11}f_{32}=0,\ \ f_{21}^2-f_{11}f_{31}=0, \ f_{22}^2-f_{12}f_{32}=0\ \  \ \ \ 
\rank
\begin{bmatrix}
 f_{11} & f_{12}\\
  f_{21} & f_{22}\\
 f_{31} & f_{32}\\
\end{bmatrix}
\leq 1,
\end{equation}
where $f_{ij}$'s are linear equations. 
In a similar way the image of the Veronese embedding $\P^2\hookrightarrow \P^5$ by degree $2$ monomials is 
given by \begin{eqnarray}
 \label{shodam2019-2}
 C& : & f_{11}f_{21}-f_{32}^2=f_{11}f_{31}-f_{22}^2=f_{21}f_{31}-f_{12}^2=0,\\ \nonumber
 & & f_{12}f_{22}-f_{31}f_{32}=f_{12}f_{32}-f_{21}f_{22}=f_{22}f_{32}-f_{11}f_{12}=0.
\end{eqnarray}
The loci of cubic fourfolds with a Veronese surface as above is $\Has_{20}$. 
%

It is now natural to generalize the quartic scroll \eqref{shodam2019} and Veronese \eqref{shodam2019-2} 
by adding $\frac{n}{2}-2$ linear equations into the defining equation of $C$ and obtain $\frac{n}{2}$ dimensional 
subvarieties of $\P^{n+1}$. We call them a quartic scroll and Veronese cycle, respectively; 
carrying in mind that they are $\frac{n}{2}$-dimensional. 
Let $\T_{\rm QS}$ (resp. $\T_{\rm V}$) be the space of smooth degree $d$ hypersurfaces 
in $\P^{n+1}$ containing a quadric scroll (resp. Veronese) cycle.  The ideal of both cycles are radical 
and written in the standard basis (Groebner basis) and so  the defining equation $f$ of a hypersurface 
containing such a cycle can be written as $f=\sum_{i=1}^{s}f_ig_i$,
where $f_i$'s are the polynomials in \eqref{shodam2019} (resp. \eqref{shodam2019-2}) and $g_i$'s are other homogeneous polynomials
such that $\deg(f_ig_i)=d$. We have computed the codimension of the image of the derivation of the maps parameterizing 
$\T_{\rm QS}$ and $\T_{\rm V}$ at many random points  for cubic hypersurfaces. These are conjecturally  the codimensions 
of $\T_{\rm QS}$ and $\T_{\rm V}$ and they are listed in Table \ref{12mar2017-2}. It would not be hard to prove that these are actual
codimensions and in the following we are going to assume this. We have used the procedures {\tt CodQuarticScroll} and
{\tt CodVeronese} for this codimension computations. In this table $L, CS, QS, V$ means Linear cycle $\P^{\frac{n}{2}}$,
cubic scroll/ruled, quadric scroll and Veronese, respectively, and the numbers below them are the codimension
of the corresponding sub loci of $\T$. The column under $M$ is the codimension of the mysterious components
of the Hodge loci that we get in Theorem \ref{hope2018}, part \ref{hope-2}. 
\begin{rem}\rm
One loci in 
Table \ref{12mar2017-2} is definitely not a component of the Hodge loci. This is namely, the loci of cubic 
six folds containing a Veronese cycle. Its codimension is $10$, whereas the upper bound for components of the
Hodge loci in this case is $8$. Assuming the Hodge conjecture in this case, it means that the Veronese cycle
is homologous to another algebraic cycle in primitive homology such that the new cycle has bigger deformation 
space. The same reasoning implies that the  loci of cubic six folds containing a Quartic scroll is a component 
of the Hodge loci as its codimension is the maximal one $h^{51}=8$. 
\end{rem}
\begin{rem}\rm
A natural generalization of $\Has_{D}$ for arbitrary cubic $n$-folds seems to have increasing codimension with respect to $D$.
If this is the case, even constructing algebraic cycles in the cubic fourfold case might not help to understand
the Hodge loci of Theorem \ref{hope2018} part \ref{hope-2} whose codimension seems to be below $\codim(\T_{\rm QS})$. 
\end{rem}
\begin{rem}\rm
For quartic scroll in \eqref{shodam2019}, by setting $f_{12}=f_{31}=0$, degenerates into a sum 
of two planes and a line 
$C_0:\P^2_1+\P^2_2+\P^1_3$:
\begin{eqnarray*}
\P^2_1 &:&  f_{21}=f_{22}^2=f_{11}=0,\\ 
\P^2_2 &:&  f_{21}^2=f_{22}=f_{32}=0,\\
\P^1_3 &:&  f_{21}^2=f_{22}^2=f_{21}f_{22}=f_{11}=f_{32}=0.
\end{eqnarray*}
Note that $\P^2_i$'s are non-reduced, the intersection $\P^2_1\cap\P^2_2$ is reduced and it is equal 
to the underlying reduced line  of $\P^1_3$. It is not clear whether the Fermat fourfold contains a quartic
scroll or its degeneration $C_0$. Note that if $C_0$ is inside a cubic fourfold $X$, 
using topological arguments it defines a class $[C_0]\in H^n_\dR(X)$, however, it is not clear
how to define it algebraically. 
\end{rem}

\begin{table}[!htbp]
\label{12mar2017-2}
\centering
\begin{tabular}{|c|c|c|c|c|c|c|c|c|}
\hline
{\tiny $\dim(X_0)$}    &  {\tiny $\dim(\T)$} & {\tiny range of codimensions} &L & CS &M &QS& V & Hodge numbers    \\ \hline
$n$     & $\bn{n+2}{3}$  & {\tiny $\binom{\frac{n}{2}+1}{3}$, $\binom{n+2}{{\rm min }\{3, \frac{n}{2}-2\}}$}  
        &&&&&   &    \tiny $h^{n,0}, h^{n-1,1},\cdots, h^{1,n-1}, h^{0,n}$ \\  \hline \hline 
 $4$    & $20$ &  $1,1$ &1&1&\fbox{1}& 1& 1  & \tiny $0,1,21,1,0$\\  \hline  
 $6$    & $56$ &  $4,8$  &4&6&\fbox{7}& 8 & \fbox{10} & \tiny $0,0,8,71,8,0,0$ \\ \hline  
  $8$   & $120$ & $10,45$  &10&16&\fbox{19}& 23 & 25 &  \tiny $0, 0, 0, 45, 253, 45, 0, 0, 0$ \\  \hline 
  $10$  & $220$ & $20, 220$ &20&32&\fbox{38}& 45 & 47  & \tiny $0, 0, 0, 1, 220, 925, 220, 1, 0, 0, 0$ \\  \hline
  $12$  & $364$ & $35,364$  &35&55&\fbox{65}& 75& 77  &   \tiny $0,0,0,0, 14,1001, 3432, 1001, 14,0,0,0,0$   \\ \hline 
\end{tabular}
\caption{Codimensions of the components of the Hodge/special loci for cubic hypersurfaces} 
\end{table}
We  started  to prepare Table \ref{16dec2018Tokyo} and Table \ref{nabeta2018} with a computer with processor 
{\tt  Intel  Core i7-7700}, $16$ GB Memory plus $16$ GB swap memory and the operating system {\tt Ubuntu 16.04}. 
It turned out that for many cases  such as $(n,m, N)=(12,3,3)$ 
in Table \ref{nabeta2018} we get the  `Memory Full' error. 
Therefore,  we had to increase the swap memory up to $170$ GB. 
Despite the low speed of the swap which slowed down the computation, the computer was 
able to use the data and give us the desired output. The computation for this example  took more than $21$ days. 
We only know that at least $18$ GB of the swap were used. 
Other time consuming computations are the cases $(n,m, N)=(10,2,4),(12,4,3)$.



 \def\cprime{$'$} \def\cprime{$'$} \def\cprime{$'$} \def\cprime{$'$}


\end{document}